\documentclass[12pt]{amsart}
\usepackage{amsmath}
\usepackage{amssymb}
\usepackage{anysize}

\marginsize{1.25cm}{1.25cm}{2.75cm}{2.75cm}

\theoremstyle{plain}
\newtheorem{theor}{Theorem}[section]
{}
\newtheorem{lemma}{Lemma}[section]

\newtheorem{corol}{Corollary}[section]

\theoremstyle{definition}
\newtheorem{defin}{Definition}[section]

\theoremstyle{remark}

\numberwithin{equation}{section}

\DeclareMathOperator{\dom}{domain}

\newcommand{\AND}{\text{ and }}
\newcommand{\IF}{\text{ if }}
 
\newcommand{\card}[1]{\lvert #1 \rvert}

\newcommand{\CH}{$2^{\aleph_0} = \aleph_1$}
\newcommand{\forces}[2]{\Vdash_{#1} \mbox{``} #2 \mbox{''}}

\newcommand{\Sacks}{{\mathbb S}}

\newcommand{\Rationals}{{\mathbb Q}}
\newcommand{\Qposet}{{\mathbb Q}}
\newcommand{\Poset}{{\mathbb P}}

\newcommand{\Naturals}{{\mathbb N}}

\newcommand{\pomega}{{\mathcal P}(\Naturals)}
\newcommand{\pomegaf}{{\mathcal P}(\Naturals)/[\Naturals]^{<\aleph_0}}
\title[Nowhere Trivial Automorphisms]{Martin's Axiom is Consistent with the Existence of
 Nowhere Trivial Automorphisms}
\author[S. Shelah]{Saharon Shelah}
\address{Department of Mathematics, Rutgers University, Hill Center,
 Piscataway, 
 New Jersey, U.S.A. 08854-8019}
\curraddr{Institute of Mathematics\\Hebrew University\\
Givat Ram, Jerusalem 91904, Israel}
\email{shelah@math.rutgers.edu}
\author[J. Stepr\={a}ns]{Juris Stepr\={a}ns}
\address{Department of Mathematics, York University,
4700 Keele Street,
Toronto, Ontario, Canada  M3J 1P3}
\curraddr{}
\email{steprans@yorku.ca}
\thanks{
The research of the first author was supported by The Israel Science
Foundation founded by the Israel Academy of Sciences and Humanities, and
by NSF grant No. NSF-DMS97-04477.
Research of the second author for this paper was partially supported by NSERC
of Canada. This is paper number 735 in the first author's personal listing}
\begin{document}
\maketitle
\begin{abstract}
Martin's~Axiom does not   imply  that all
automorphisms of $\pomegaf$ are somewhere trivial. An alternate method for
obtaining models where every automorphism of $\pomegaf$ is somewhere trivial
is explained.
\end{abstract}
\bibliographystyle{plain}
\section{Introduction}
In \cite{veli.oca} Veli{\v{c}}kovi{\'c} constructed a model of Martin's~Axiom
in which  there is a non-trivial automorphism of $\pomegaf$. As well as
answering a question posed in \cite{veli.def}, this put into
context
another result of the \cite{veli.oca}  showing that the conjunction of MA and
OCA implies that all 
automorphisms are trivial.
However, the non-trivial automorphims constructed by Veli{\v{c}}kovi{\'c} is
trivial on many infinite subsets of the integers. Indeed, it was shown in
\cite{step.28} that this is unavoidable since
 every automorphism of $\pomegaf$ is somewhere
trivial in Veli{\v{c}}kovi{\'c}'s model of \cite{veli.oca}.

Hence, the question arises of whether or not Martin's~Axiom alone is
sufficient to   imply  that, while  there may be non-trivial
automorphisms, nevertheless, all
automorphisms of $\pomegaf$ are somewhere trivial.
The main result of this paper is that this is not the case.

The last section presents a simple, alternate method for obtaining models
where all automorphisms are somewhere trivial. It has the advantage
that it can produce models where ${\mathfrak d} = \aleph_1$ whereas the 
oracle chain condition method adds Cohen reals and so does not achieve this. 

\section{Martin's Axiom and a nowhere trivial automorphism}
If $\alpha$ and $\beta $ are ordinals then the notation $[\alpha,
\beta)$ will be used to denote the set $\beta \setminus \alpha$.
The relations $\equiv^*$, $\subseteq^*$ and $\supseteq^*$ will have the usual
meaning as relations on subsets of the integers modulo a finite set. The
convention on forcing partial orders will be that larger conditions force more
information. 
\begin{defin}\label{d:p}
If $W$ is a set of ordinals then the indexed family  ${\mathfrak S} = \{(A_\xi, F_\xi,
{\mathfrak B}_\xi)\}_{\xi \in W}$ 
will be said to be a {\em tower of permutations} if 
\begin{enumerate}
\item $A_\xi \subseteq \Naturals$ and $F_\xi$ is a
permutation of $\Naturals$  for each $\xi$
\item  $F_\xi\restriction m$ is a
permutation of $m$  for each $m$ in $A_\xi$
\item ${\mathfrak B}_\xi$ is a finite
subalgebra of ${\mathcal P}(\Naturals)$  for each $\xi$
\item if $\xi \in \zeta$ then ${\mathfrak B}_\xi \subseteq {\mathfrak
B}_\zeta$ and $A_\xi \supseteq ^* A_\zeta$ 
\item if $\xi \in \zeta$ then $F_\zeta(B) \equiv^* F_\xi(B)$ for each 
$B \in {\mathfrak B}_\xi$.
\end{enumerate}
Define $\Qposet({\mathfrak S})$
to be the set consisting of all quadruples $p = (a^p,
f^p,\alpha^p,  {\mathfrak B}^p) $ such that 
\begin{enumerate}
\item $a^p \subseteq \Naturals$ is a finite subset 
\item $f^p\restriction m$ is a permutation of $m$ for each  $m \in a^p$
and the domain of $f^p$ is $\max(a^p)$
\item $\alpha^p \in W$ 
\item \label{i:9} $\max(a^p) \in A_{\alpha^p}$ 
\end{enumerate}
and the relation $\leq$ on $\Qposet({\mathfrak S})$ is defined by $p \leq q$ if and only if
\begin{equation}\label{e:d1}
a^p \subseteq a^q,\ \ 
 (\max(a^p)+1)\cap a^q = a^p,\ \ 
f^p \subseteq f^q,\ \ 
{\mathfrak B}^p \subseteq {\mathfrak B}^q,\ \ 
\alpha^p \subseteq \alpha^q
\end{equation}
\begin{equation}\label{e:d5}
(A_{\alpha^q} \setminus \max(a^q))  \cup (a^q \setminus a^p )\subseteq A_{\alpha^p}
\end{equation}
and, for each $B$ belonging to  ${\mathfrak B}^p\cap{\mathfrak B}_{\alpha^p}$,
 the following two conditions hold:
\begin{equation}\label{e:d7}
(\forall \{n,m\} \in [a^q\setminus \max(a^p)]^2)
f^q(B \cap [m,n)) = 
F_{\alpha^p}(B \cap [m, n))
\end{equation}
\begin{equation}\label{e:d8}
(\forall  \{n,m\} \in [A_{\alpha^q}\setminus
\max(a^q)]^2) F_{\alpha^q}(B \cap [m,n)) =
F_{\alpha^p}(B \cap [m, n)) .
\end{equation}
If $G$  
is generic for $\Qposet( {\mathfrak S})$
then define $A_{{\mathfrak S}}[G] = \bigcup_{p\in G}a^p$ 
 and  
$F_{{\mathfrak S}}[G]$   to be 
$\bigcup_{p \in G}f^p$.
\end{defin}

\begin{lemma}\label{l:tr}
For any tower of permutations $\mathfrak S$ the structure
$(\Qposet({\mathfrak S}), \leq)$ is a partial order.
\end{lemma}
\begin{proof} 
That $(\Qposet({\mathfrak S}), \leq)$ is reflexive and antisymmetric
is obvious. To prove transitivity suppose that $p \leq q$ and $q \leq
r$.
The condition~\ref{e:d1} for $p \leq r$ is easily seen to be satisfied.
To see that condition~\ref{e:d5} for $p \leq r$ is  satisfied note that
$$A_{\alpha^r} \setminus \max(a^r) \subseteq A_{\alpha^q} \setminus
\max(a^r) \subseteq A_{\alpha^q} \setminus \max(a^q) \subseteq
A_{\alpha_p}$$
and that 
$$a^r \setminus a^p \subseteq (a^r \setminus a^q) \cup (a^q \setminus
a^p) \subseteq A_{\alpha^q} \setminus \max(a^q) \cup A_{\alpha^p}
\subseteq A_{\alpha^p}  $$
which shows that $(A_{\alpha^r} \setminus \max(a^r))\cup (a^r
\setminus a^p) \subseteq A_{\alpha^p}$, as required.

To show that conditions~\ref{e:d7} and~\ref{e:d8} hold, 
let $B \in
{\mathfrak B}^p\cap {\mathfrak B}_{\alpha^p}$. 
Given any pair $\{n,m\} \in [a^r\setminus \max(a^p)]^2$,
it may, without loss of generality be assumed that 
$n$ and $m$ are successive elements of  $a^r\setminus \max(a^p)$.
Hence, either $\{n,m\} \in [(a^r\setminus \max(a^q)]^2$ or
$\{n,m\} \in [(a^q\setminus \max(a^p)]^2$. In the second case, it
follows immediately from the fact that  $p \leq q$ that 
$f^q(B\cap[n,m)) 
= F_{\alpha^p} (B\cap[n,m)) $. 
Since $q \leq r$ it follows that $f^q
\subseteq f^r$ and so $f^r(B\cap[n,m)) 
= F_{\alpha^p} (B\cap[n,m)) $. On the other hand, in the first case
$f^r(B\cap[n,m)) 
= F_{\alpha^q} (B\cap[n,m)) $ since $q \leq r$ and ${\mathfrak B}_{\alpha^p}
\subseteq {\mathfrak B}_{\alpha^q}$ . Moreover, using
condition~\ref{e:d8} and $p \leq q$ it is possible to conclude that 
 $F_{\alpha^q} (B\cap[n,m)) 
= F_{\alpha^p} (B\cap[n,m)) $. Hence, in either case
$f^r(B\cap[n,m)) 
= F_{\alpha^p} (B\cap[n,m)) $ which establishes that
condition~\ref{e:d7} holds for $p \leq r$.

To see  that
condition~\ref{e:d8} holds for $p \leq r$ 
let  $\{n,m\} \in [A_{\alpha^r}\setminus \max(a^r)]^2$. It follows
from $p \leq q$ and $q\leq r$
that $F_{\alpha^q} (B\cap[n,m)) 
= F_{\alpha^r} (B\cap[n,m)) $ and, since $\{n,m\} \subseteq
A_{\alpha^r}\setminus \max(a^r) \subseteq A_{\alpha^q}\setminus \max(a^q)$,
  that
$F_{\alpha^q }(B\cap[n,m)) 
= F_{\alpha^p} (B\cap[n,m)) $. Hence, 
$F_{\alpha^r} (B\cap[n,m)) 
= F_{\alpha^p} (B\cap[n,m)) $ establishing condition~\ref{e:d8}.
\end{proof}

\begin{lemma}\label{l:density}
Given a tower of permutations ${\mathfrak S} = \{(A_\xi, F_\xi,
{\mathfrak B}_\xi)\}_{\xi \in W}$,
an integer $n$, $B\in \bigcup_{\xi \in W}{\mathfrak B}_{\xi}$ and  $\zeta \in W$
the following
sets are dense in $ \Qposet({\mathfrak S})$:
\begin{equation}\label{e:ds1}
\{p \in \Qposet({\mathfrak S}) :
 \max(a^p) > n\}
\end{equation}
\begin{equation}\label{e:ds2}
\{p \in \Qposet({\mathfrak S}) :
 \alpha^p \geq \zeta\}
\end{equation}
\begin{equation}\label{e:ds3}
\{p \in \Qposet({\mathfrak S}) :
 B \in {\mathfrak B}^p \}
\end{equation}
\end{lemma}
\begin{proof}
To prove that the set~\ref{e:ds1} is dense let $p \in
\Qposet({\mathfrak S})$ and $n$ be given. Let $k \in A^{\alpha^p}$ be
such that $k > n$. Using Condition~\ref{i:9} of Definition~\ref{d:p}
it follows that $\max(a^p) \in A_{\alpha^p}$ and, hence,
$F_{\alpha_p}\restriction [\max(a^p),k)$ is a permutation of
$[\max(a^p),k)$. Letting $q = (a^p\cup\{k\}, f^p \cup
F_{\alpha_p}\restriction [\max(a^p),k), \alpha^p, {\mathfrak B}^p)$ it
follows that $q \geq p$.
Observe for later reference, that is has actually been shown that
\begin{equation}\label{l:dexx}
(\forall p \in \Qposet({\mathfrak S}))(\forall n \in \Naturals)(\exists q \geq
p) max(a^q) > n \AND \alpha^q = \alpha^p \AND {\mathfrak B}^p = {\mathfrak B}^q.
\end{equation}

To prove that the set~\ref{e:ds2} is dense let $p \in
\Qposet({\mathfrak S})$ and $\zeta\in W$ be given. Since it may as well be
assumed that $\zeta > \alpha^p$, it is possible to find  $n$ 
so large that $A_\zeta \setminus n \subseteq A_{\alpha^p}$ and
for all $\{i,j\} \in [A_\zeta\setminus n]^2$ and $B\in {\mathfrak B}^p$
$$F_\zeta(B\cap [i,j)) =  F_{\alpha^p}(B\cap [i,j)).$$
Using the set~\ref{e:ds1} of Lemma~\ref{l:density} choose $q$ such that $ p \leq q$
and $n < \max(a^q)$. From \ref{l:dexx} it can be assumed that ${\mathfrak
B}^q = {\mathfrak B}^p$ and that $\alpha^q = \alpha^p$.
Now, let $r = (a^q, f^q, \zeta,
{\mathfrak B}^p)$. That conditions~\ref{e:d5} and \ref{e:d8}
for the relation  $q \leq r$ is
satisfied follows from the choice of $n$ while condition~\ref{e:d1} is
obvious.  Condition~\ref{e:d7} has no
content in the case of $q \leq r$ since $a^q = a^r$. Now use transitivity and
the fact that 
$p \leq q$. Observe for later reference, that is has actually been shown that
\begin{equation}\label{l:dex}
(\forall p \in \Qposet({\mathfrak S}))(\forall \zeta > \alpha^p)(\exists q \geq
p) \alpha^q = \zeta .
\end{equation}

There is no problem in proving that the set~\ref{e:ds3} is dense. 
\end{proof}

\begin{lemma}\label{l:genericity1}
If ${\mathfrak S} = \{(A_\xi, F_\xi, {\mathfrak B}_\xi)\}_{\xi \in W}$ is a tower of permutations 
and $p \in \Qposet({\mathfrak S})$ then
\begin{equation}\label{e:g1}
p\forces{\Qposet({\mathfrak S})}{A_{{\mathfrak S}}[G] \setminus \max(a^p)\subseteq
A_{\alpha^p}}
\end{equation}
\begin{equation}\label{e:g2}p\forces{\Qposet({\mathfrak S})}{
(\forall B \in {\mathfrak B}^p\cap {\mathfrak B}_{\alpha^p})(\forall \{n,m\} \in [A_{\alpha^p}\setminus \max(a^p)]^2) 
F_{{\mathfrak S}}[G](B\cap [n,m)) = F_{\alpha^p}(B\cap [n,m))}
\end{equation}
\end{lemma}
\begin{proof}
This is standard using condition~\ref{e:d5} for \ref{e:g1} and condition~\ref{e:d7} for \ref{e:g2}.
\end{proof}

Let $\kappa$ be a regular uncountable cardinal and let $C\subseteq
\kappa$ be any set containing 0 and 
closed under limits of increasing $\omega_1$-sequences such that $\kappa \setminus C$ is
also unbounded.  Now define $\Poset_\eta$, as well as a
$\Poset_\eta$-name for a tower of permutations 
${\mathfrak S}_\eta = \{(A_\zeta, F_\zeta,{\mathfrak B}_\zeta)\}_{\zeta\in
  \eta\cap C}$,  
 by induction on $\eta$. 
Let $ A_0 \in [\Naturals]^{\aleph_0}$
 and $F_0$ be arbitrary subject to the fact that 
$F_0$ is a permutation of  $\Naturals$ 
such that $F_0\restriction [n,m)$ is a permutation for each $\{n,m\}\subseteq
A_0$. 
Let  ${\mathfrak B}_0 = {\mathcal P}( \Naturals)$ in the sense of
the ground model. Then  let 
${\mathfrak S}_{1} = \{(A_0, F_0, {\mathfrak B}_0)\}$
 and let $\Poset_1 = \Qposet_0$ be Cohen forcing. 
If $\eta$ is a limit then $\Poset_\eta$ is simply the finite support
limit of $\{\Poset_\zeta\}_{\zeta\in\eta}$ and ${\mathfrak S}_\eta
= \bigcup_{\xi \in C \cap
\eta}{\mathfrak S}_\xi$.
 If $\eta \notin C$ then  
$\Poset_{\eta + 1} = \Poset_\eta*\Qposet_\eta$ where $\Qposet_\eta$ is
a ccc partial  order
chosen according to some bookkeeping scheme which will guarantee that
Martin's Axiom holds at stage $\kappa$. In this case
${\mathfrak S}_{\eta + 1} = {\mathfrak S}_{\eta}$.
 If $\eta \in C\setminus \{0\}$ then  
$\Poset_{\eta + 1} = \Poset_\eta*\Qposet({\mathfrak S}_{\eta})$ and $A_\eta$ is defined  to be 
$A_{{\mathfrak S}_{\eta}}[G]$,  
$F_\eta$ is defined  to be 
$F_{{\mathfrak S}_{\eta}}[G]$ where $G$ is the
canonical name for the generic set on $\Qposet({\mathfrak
S}_{\eta})$. In this case
${\mathfrak S}_{\eta + 1} = {\mathfrak S}_{\eta}\cup \{(A_\eta,
F_\eta,{\mathcal P}(\Naturals)\cap V^{\Poset_{\eta}})\}$. As usual, if $p \in
\Poset_\eta$ then $p\restriction \alpha\forces{\Poset_\alpha}{p(\alpha) \in
  \Qposet_\alpha}$. 

\begin{defin}\label{d:prime}
Let $\Poset_{\alpha}^w \subseteq
\Poset_\alpha$ consist of all those $p \in \Poset_\alpha$  such that  
there are $k \in \Naturals$, $\Gamma \in [C\cap \alpha]^{<\aleph_0}$ and
$\{(a_\gamma, f_\gamma)\}_{\gamma\in\Gamma}$
such that 
\begin{itemize}
\item $\Gamma = \dom(p) \cap C \setminus \{0\}$
\item $0\in\Gamma$
\item if  $\gamma\in\Gamma$ 
then $p \restriction \gamma\forces{\Poset_\gamma}{p(\gamma) =
(\check{a}_\gamma, \check{f}_\gamma, \check{\alpha}_\gamma, {\mathfrak
  B}^\gamma)}$ for some 
${\mathfrak B}^\gamma$
\item $\alpha_\gamma \in\Gamma\cap
\gamma$ for each $\gamma\in\Gamma$
\item if  $\gamma\in\Gamma$ then $\max(a_\gamma) = k$
\item if $\gamma$ and $\gamma'$ are in $ \Gamma$ and $\gamma' <
\gamma$ then $p\restriction \gamma\forces{\Poset_{\gamma}}{{\mathfrak B}^{\gamma'} =
  {\mathfrak B}^{\gamma}\cap  {\mathfrak B}_{\gamma'} }$. 
\end{itemize}
The pair $(k, \{(a_\gamma, f_\gamma)\}_{\gamma \in \Gamma})$
will
be said to witness that $p \in 
\Poset_{\alpha}^w$. 
Let $\Poset_{\alpha}^* \subseteq
\Poset_\alpha$ consist of all those $p \in \Poset_\alpha^w$  such that,  
in addition to the other requirements,  $\alpha_\gamma = \max(\Gamma\cap
\gamma)$ for each $\gamma\in\Gamma$.
\end{defin}
\begin{lemma}\label{l:04}
If $p \in \Poset_\alpha^w$ and this is witnessed by $(k, \{(a_\gamma,
f_\gamma)\}_{\gamma \in \Gamma})$ then
\begin{itemize}
\item if $\gamma\in\Gamma$   then
$p\restriction \gamma + 1\forces{\Poset_{\gamma + 1}}{A_{\alpha_\gamma} \supseteq
(A_{\gamma} \setminus k)}$
\item if $\gamma\in\Gamma$   then
$p\restriction \gamma + 1\forces{\Poset_{\gamma + 1}}{(\forall B
\in {\mathfrak B}^\gamma\cap {\mathfrak B}_{\alpha_\gamma})(\forall  \{n,m\} 
\in [A_\gamma \setminus k]^2)
F_{\gamma}(B \cap [n,m)) = F_{\alpha_\gamma}(B\cap [n,m))}$
\end{itemize}
\end{lemma}
\begin{proof}
This is an immediate consequence of Lemma~\ref{l:genericity1}.  
\end{proof}

\begin{defin}
  If $p \in \Poset_\alpha^w$ is witnessed by  $(k, \{(a_\gamma,
  f_\gamma)\}_{\gamma \in \Gamma})$ then define $p^+\in \Poset_\alpha^*$ 
by
$$p^+(\xi) = \begin{cases}
p(\xi) & \IF \xi \notin \Gamma\\
(a_\xi, f_\xi, \max(\Gamma\cap \xi), {\mathfrak B}^{p(\xi)}) & \IF \xi \in
\Gamma .
\end{cases}$$
\end{defin}

\begin{lemma}\label{l:o6}
If $p \in \Poset_\alpha^w$ then $p^+ \geq p$.  
\end{lemma}
\begin{proof}
Proceed by induction on $\beta \in \alpha$ to show that $p^+\restriction \beta
  \geq p\restriction \beta$. Note that the cases $\beta = 0$ or $\beta$ a
  limit pose no problem. Given that $p^+\restriction \beta
  \geq p\restriction \beta$ let $\Gamma \cap \beta$ be enumerated, in order,
by $\{\gamma_1, \gamma_2,\ldots \gamma_n\}$.
  If follows directly from Lemma~\ref{l:04} and the definition of $p^+$
that
\begin{itemize}
\item 
$p^+\restriction \beta\forces{\Poset_{\beta}}{A_{\gamma_1}\setminus \check{k} \supseteq
(A_{\gamma_2}\setminus \check{k})\ldots \supseteq
(A_{\gamma_n}\setminus \check{k})}$
\item 
$p^+\restriction \beta\forces{\Poset_{\beta}}{(\forall B
\in {\mathfrak B}^{\gamma_j}\cap {\mathfrak B}_{\gamma_i})(\forall  \{m,m'\}
\in A_{\gamma_j} \setminus \check{k}) 
F_{\gamma_j}(B \cap [m,m')) = F_{\gamma_i}(B\cap [m,m'))}$ if $i \leq j$.
\end{itemize}
In particular, noting that there is some $i$ such that
 $\alpha_\beta = \gamma_i$,
$$p^+\restriction \beta\forces{\Poset_{\beta}}{A_{\alpha_\beta}\setminus \check{k} \supseteq
A_{\gamma_n}}$$ and
$$p^+\restriction \beta\forces{\Poset_{\beta}}{(\forall B
\in {\mathfrak B}^{\alpha_\beta}\cap {\mathfrak B}_{\gamma_n})(\forall
\{m,m'\} \in A_{\alpha_\beta} \setminus k) 
F_{\alpha_\beta}(B \cap [m,m')) = F_{\gamma_n}(B\cap [m,m'))} .$$
Hence, $p^+\restriction \beta + 1
  \geq p\restriction \beta+1$.
\end{proof}


\begin{lemma}\label{l:nice}
For each $\alpha \leq \kappa$ the subset $\Poset_{\alpha}^w$ is dense in
$\Poset_\alpha$.
\end{lemma}
\begin{proof}
Proceed by induction on $\alpha$ noting that the cases $\alpha \leq 1$
and $\alpha$ a limit are trivial. Therefore, suppose that
Lemma~\ref{l:nice} has been established for $\beta$ and that $p \in
\Poset_{\beta+1}$.
Without loss of generality it may be assumed that $\beta
\in C$. Choose $p' \geq p\restriction \beta$ and
$(a, f, \zeta)$ such that $\zeta \in \beta$ and 
 $p'\forces{\Poset_\beta}{p(\beta) = (\check{a},
\check{f}, \check{\zeta}, {\mathfrak B})}$ and, moreover, 
 $p'\forces{\Poset_\beta}{{\mathfrak B} = \{B_i\}_{i\in m}}$ and, for
each $i\in m$ the condition $p'$ decides the value of the least
ordinal $\zeta(i) \in \beta$ such that
$B_i \in {\mathfrak B}_{\zeta(i)}$. Furthermore,
it may, without loss of generality be assumed that $\{\zeta\}
\cup\{\zeta(i)\}_{i\in m} 
\subseteq \dom(p')$ and that $p'\restriction \zeta(i)
\forces{\Poset_{\zeta(i)}}{B_i \in {\mathfrak B}^{p'(\zeta(i))}}$ for $i \in
m$. 
Then use the induction
hypothesis to find $p'' \in \Poset_{\beta}^w$ extending  $p'$ and such that 
 $(k, \{(a_\gamma, f_\gamma)\}_{\gamma \in \Gamma})$
witnesses 
that $p'' \in \Poset_{\beta}^w$. Without loss of generality, $k \geq max(a)$.
Now, the fact that $\{\max(a), k\} \subseteq  A_\zeta$ 
implies that $F_\zeta\restriction
[\max(a), k)$ is a permutation of $[\max(a), k)$. Hence it is possible
to define $f' = f \cup F_\zeta\restriction
[\max(a), k)$ and $a' = a \cup \{k\}$. Then define $q$ so that
$$q(\xi) = \begin{cases}
p''(\xi) & \IF \xi \in \beta\\
(a',f',\zeta,{\mathfrak B}\cup \bigcup_{\gamma \in \Gamma}{\mathfrak B}^{p'(\gamma)}) & \IF \xi = \beta 
           \end{cases}$$
and note that $q \in \Poset^w_{\beta + 1}$ and $q\geq p$.
\end{proof}

\begin{corol}
  For each $\alpha \leq \kappa$ the subset $\Poset_{\alpha}^*$ is dense in
$\Poset_\alpha$.
\end{corol}
\begin{proof}
  This is immediate from Lemma~\ref{l:nice} and Lemma~\ref{l:o6}.
\end{proof}

\begin{lemma}\label{l:add}
Given that  $p \in \Poset_{\alpha}^w$ and that this is witnessed by
$(k, \{(a_\gamma,f_\gamma)\}_{\gamma\in \Gamma})$ and $\mu
\in (C\cap \alpha) \setminus \dom(p)$ then the following condition
$p\langle\mu\rangle$ extends $p$ and is also in  
$\Poset_{\alpha}^w$:
$$p\langle\mu\rangle(\xi) \begin{cases}
p(\xi) & \IF \xi \neq \mu\\
(\{k\},\{(j,j)\}_{j\in k}, 0, \bigcup_{\gamma\in
\Gamma\cap \mu}{\mathfrak B}^{p(\gamma)}) & \IF \xi = \mu .
\end{cases}$$
\end{lemma}
\begin{proof}
Notice that since $\mu \notin \dom(p)$ the restrictions on extension do
not apply and it is easy to check that the condition
$(\{k\},\{(j,j)\}_{j\in k}, \max(\Gamma\cap \mu), \bigcup_{\gamma\in 
\Gamma\cap \mu}{\mathfrak B}^{p(\gamma)})$ belongs to $\Qposet_\mu$.
 Since $p \in \Poset^w_\alpha$ it follows that $k \in A_0$ and, hence, 
that  $p\langle\mu\rangle \in \Poset_\alpha^w$. That
$p\langle\mu\rangle \geq p$ is immediate from the definition.
\end{proof}

\begin{lemma}\label{l:comb}
Suppose that $p$ and $q$ are  conditions in $\Poset_\alpha$
such that:
\begin{itemize}
\item   $p\in \Poset_{\alpha}^*$
is witnessed by $(k, \{(a_p^\gamma, f_p^\gamma)\}_{\gamma \in \Gamma_p})$
\item  $q \in \Poset_{\alpha}^*$
is witnessed by $(k, \{(a_q^\gamma, f_q^\gamma)\}_{\gamma \in \Gamma_q})$
\item $\max(\dom(p)) = \max(\dom(q)) = \bar{\gamma} \in \Gamma_p \cap \Gamma_q$
\item $(a_q^{\bar{\gamma}},f_q^{\bar{\gamma}}) = (a_p^{\bar{\gamma}},f_p^{\bar{\gamma}})$
\item $\max(\dom(q))\cap \bar{\gamma} < \min(\dom(p)\setminus \{0\})$
\item $p(0)$ and $(0)$ are compatible.
\end{itemize}
Under these conditions $p$ and $q$ are compatible.
\end{lemma}
\begin{proof}
%
Let $\sigma$ the maximum member of $\Gamma_{{q}}\setminus
\{\bar{\gamma}\}$.
Define $q \sqcup p$ by
$$(q\sqcup p)(\xi) = \begin{cases}
q(\xi) & \IF \xi \in \dom(q) \setminus (\Gamma_q\cup \mu)\\
p(\xi) & \IF \xi \in \dom(p) \setminus (\Gamma_p\cup \mu)\\
(a ^{{q}(\xi)},f^{{q}(\xi)},\alpha^{{q}(\xi)},({\mathfrak
  B}^{{p}}\cap {\mathfrak B}_0)\cup {\mathfrak
  B}^{{q}(\xi)}) & \IF \xi \in \Gamma_q \\
(a^{{p}(\xi)},f^{{p}(\xi)},\alpha^{{p}(\xi)},{\mathfrak
  B}^{{p}(\xi)}\cup  
{\mathfrak  B}^{{q}(\sigma)}) & \IF \xi \in \Gamma_p \\
(a^{{p}(\bar{\gamma})}, f^{{p}(\bar{\gamma})}, \alpha^q(\bar{\gamma}),
{\mathfrak B}^{{p}(\bar{\gamma})}  \cup {\mathfrak B}^{q(\bar{\gamma})}) &
\IF \xi = \bar{\gamma}
                   \end{cases}$$
and define $p \sqcup q$ by
$$(p\sqcup q)(\xi) = \begin{cases}
q(\xi) & \IF \xi \in \dom(q) \setminus (\Gamma_q\cup \mu)\\
p(\xi) & \IF \xi \in \dom(p) \setminus (\Gamma_p\cup \mu)\\
(a ^{{q}(\xi)},f^{{q}(\xi)},\alpha^{{q}(\xi)},({\mathfrak
  B}^{{p}}\cap {\mathfrak B}_0)\cup {\mathfrak
  B}^{{q}(\xi)}) & \IF \xi \in \Gamma_q \\
(a^{{p}(\xi)},f^{{p}(\xi)},\alpha^{{p}(\xi)},{\mathfrak
  B}^{{p}(\xi)}\cup  
{\mathfrak  B}^{{q}(\sigma)}) & \IF \xi \in \Gamma_p \\
(a^{{p}(\bar{\gamma})}, f^{{p}(\bar{\gamma})}, \alpha^p(\bar{\gamma}),
{\mathfrak B}^{{p}(\bar{\gamma})}  \cup {\mathfrak B}^{q(\bar{\gamma})}) &
\IF \xi = \bar{\gamma}
                   \end{cases}$$
(The only difference is to be found in the last lines of the two definitions.)
It is easy to check that
 $p \leq p \sqcup q$ and $q\leq q \sqcup p$ and that both
 $ p \sqcup q$ and $ q \sqcup p$ belong to $\Poset_\alpha^w$.
Moreover  $( p \sqcup q)^+ = ( q \sqcup p)^+$. 
 Hence $q$ and $p$ are compatible.
\end{proof}
 
\begin{lemma}\label{l:esay}
If the cofinality of $\zeta\cap C$ is not $\omega_1$ and
${\mathfrak S} = \{(A_\xi,F_\xi,{\mathfrak B}_\xi)\}_{\xi\in \zeta
\cap C} $ is a tower of permutations then $\Qposet({\mathfrak S})$ has
property K.
\end{lemma}
\begin{proof}
This is standard. 
If the cofinality of $\zeta\cap C$ is less than $\omega_1$ then
choose a countable, cofinal subset $C'$ of  $\zeta\cap C$. 
Using observation~\ref{l:dex} of Lemma~\ref{l:density} it follows that
the set of all $p \in \Qposet({\mathfrak S})$ such that $\alpha^p \in
C'$ is dense. It follows that $\Qposet({\mathfrak S})$ has a
$\sigma$-centred dense subset.

On the other hand, if the cofinality of $\zeta\cap C$ is greater than
$\omega_1$ 
and  $\{p_\xi\}_{\xi\in\omega_1} \subseteq
\Qposet({\mathfrak S})$ it is possible to choose $\theta \in C\cap
\zeta$ such that $\alpha^{p_\xi} < \theta$ for each $\xi$.
Using observation~\ref{l:dex} of Lemma~\ref{l:density} it may be
assumed, by extending each condition, that $\alpha^{p_\xi} = \theta$ for each $\xi$.
 Now there are $a$ and
$f$ and an uncountable set of $\xi$ such that  
 $a^{p_\xi} = a$ and $f^{p_\xi} = f$. 
Any two of these are easily seen to be compatible.
\end{proof}

The fact that the tower of permutations needs to be generic, or at least
some other condition must be satisfied, in order  for
Lemma~\ref{l:esay} to hold has been observed in Theorem~2 of \cite{step.15}.

\begin{lemma}\label{l:2}
$\Poset_\alpha$ has the countable chain condition for each $\alpha$.
\end{lemma}
\begin{proof}
Proceed by induction on $\alpha$. If $\alpha = 1$ the result is
immediate and if $\alpha$ is a limit the result follows from
the induction hypothesis and the finite support of the iteration.
Therefore consider the case $\alpha = \beta + 1$ and assume that the
countable chain condition has already been established for
$\Poset_\beta$.

Next, observe that if $C\cap \beta$ is not cofinal in $\beta$ then the
induction hypothesis is easily applied since, in this case, $C$ has a maximal
element below $\beta$ and so $\Qposet_\beta$ has the countable chain condition
 by Lemma~\ref{l:esay}.
If $C\cap \beta$ has cofinality different from
$\omega_1$ then, once again Lemma~\ref{l:esay} implies that
$\Qposet_\beta$ has the countable chain condition; therefore, in either case,
 so does
$\Poset_{\beta+1}$. Hence it remains to consider the case that 
$C\cap \beta$ is cofinal in $\beta$ and has cofinality $\omega_1$. From the hypothesis on $C$,
and the fact that $\beta$ must be the limit of $C\cap \beta$ it
follows that $\beta \in C$. 
By appealing to
Lemma~\ref{l:nice} and extending the conditions in question,  it is
possible to guarantee that each $p_\xi$ is in $\Poset_{\alpha}^*$ and that
this is witnessed by  $(k_\xi, \{(a_\xi^\gamma, f_\xi^\gamma)\}_{\gamma \in
  \Gamma_\xi})$
 As well, by thinning out, it may be
assumed that there is $k$ such that $k_\xi = k$ for each $\xi$, and there is a
pair $(a,f)$ such that $(a_\xi^\beta, f_\xi^\beta) = (a,f)$ for each $\xi$
and that
$\{\dom(p_\xi)\}_{\xi \in \omega_1}$ form a $\Delta$-system with root
$\Delta$. 
Using the fact that $\beta$ is a limit, choose some $\mu\in \beta\cap C$ such
that there is some uncountable $\Lambda \subseteq \omega_1$ such that
if $\{\eta, \zeta\} \in [\Lambda]^2$ and  $\eta \in \zeta$ then 
$\max(\dom(p_\eta \cap \beta)) < \min(\dom(p_\zeta \cap [\mu,\beta)))$.
It may as well be assumed that $\mu \notin \dom(p_\eta)$ for each $\eta \in
\Lambda$. 
Then use Lemma~\ref{l:add} to extend each $p_\eta$ to some
$p_\eta\langle\mu\rangle^+$. 

Now, using the induction hypothesis, find $\{\eta, \zeta\} \in
[\Lambda]^2$ such that  $\eta \in \zeta$ and $p_\eta\langle\mu\rangle^+\restriction\mu + 1$ is
compatible with $p_\zeta\langle\mu\rangle^+\restriction\mu + 1$. Choose $r \in
\Poset_{\mu + 1}$ 
extending both $p_\eta\langle\mu\rangle^+\restriction\mu + 1$ and 
$p_\zeta\langle\mu\rangle^+\restriction\mu + 1$.

Let $G\subseteq \Poset_{\mu+1}$ be generic over $V$ and containing $r$.
Observe that $\Poset_\alpha/\Poset_{\mu + 1}$ as interpreted in $V[G]$ is a
partial order like $\Poset_\alpha$ in $V$ except that the first factor of
$\Poset_\alpha/\Poset_{\mu + 1}$ is $\Qposet(\{(A_\mu,
F_\mu,\pomega\cap V[G])\})$. 
Since $r\forces{\Poset_{\mu + 1}}{k \in A_\mu}$  it follows
that $k \in A_\mu$, in $V[G]$.
It follows that
 $(p_\eta\langle\mu\rangle^+\restriction [\mu,\alpha))$ and
 $(p_\zeta\langle\mu\rangle^+ \restriction [\mu,\alpha))$ belong to 
$(\Poset_\alpha/\Poset_{\mu + 1})^*$. Now use Lemma~\ref{l:comb} in $V[G]$
to conclude that 
 $(p_\eta\langle\mu\rangle^+ \restriction [\mu,\alpha))$
and  $(p_\zeta\langle\mu\rangle^+ \restriction [\mu,\alpha))$ are compatible.
Hence, so are
 $p_\eta$ and $p_\zeta$.
\end{proof}

\begin{lemma}\label{l:kl}
Let ${\mathfrak S} = \{(A_\xi, F_\xi,
{\mathfrak B}_\xi)\}_{\xi \in W}$ be a tower of permutations and
suppose that $A \in [\Naturals]^{\aleph_0}$ and $\psi$ is a one-to-one
function from  $A$ to $\Naturals$. 
Then $\Qposet({\mathfrak S})$ forces that there are infinitely many
$a\in A$ such that  $F_{\mathfrak S}[G](a) \neq \psi(a)$.
\end{lemma}
\begin{proof}
Let $p\in \Qposet({\mathfrak S})$ and suppose that
$p\forces{\Qposet({\mathfrak S})}{
(\forall a \in \check{A}\setminus  \check{k})\psi(a) = F_{\mathfrak
S}[G](a)}$ for some integer $k$. From Lemma~\ref{l:density} it can be
assumed that $k \leq \max(a^p)$.
Now choose  $m \in A_{\alpha^p}$ so large that there exist distinct
integers $i$ and $j$ in $[\max(a^p),m)\cap A$ such that
 $i \in B$ if and only if $j \in B$ for all $B \in
{\mathfrak B}^p$. It is possible to choose a bijection $g : \{i,j\} \to
\{F_{\alpha^p}(i), F_{\alpha^p}(j)\}$ such that $g \neq \psi
\restriction \{i,j\}$. However, letting $f = F_{\alpha^p}\restriction
([\max(a^p),m) \setminus \{i,j\}) \cup g$ it follows that $f(B \cap
[\max(a^p),m)) = F_{\alpha^p}(B \cap [\max(a^p),m))$ for each $B \in
{\mathfrak B}^p$. Hence $p \leq q = (a^p\cup \{m\}, f^p \cup f, \alpha^p,
{\mathfrak B}^p)$ and $q\forces{\Qposet({\mathfrak S})}{
\check{\psi}(\check{i}) \neq F_{\mathfrak
S}[G](\check{i})} $.
\end{proof}

 Lemma~\ref{l:esay} and Lemma~\ref{l:kl}, in the case
 case when $C\cap \zeta$ has the
maximal element, together imply that there is  a trivial
automorphism $\Phi$ of $\pomegaf$ and a $\sigma$-centred forcing 
$\Poset$ such that $\Phi$ can be extended to a
trivial  automorphism of $\pomegaf$ in two different ways in the generic
 extension by $\Poset$. Hence, these lemmas can be thought of as 
strengthening
 the folklore result that certain automorphisms --- such as the identity ---
 can be extended to generic trivial automorphisms by the natural
 $\sigma$-centred  forcing; in other words, the forcing which approximates the
 generic permutation by finite permutations and promises to respect the
 identity on certain, finitely many,  infinite sets.  

\begin{lemma}\label{l:kl+}
Let ${\mathfrak S} = \{(A_\xi, F_\xi,
{\mathfrak B}_\xi)\}_{\xi \in W}$ be a tower of permutations such that $W$ has
a maximal element $\theta$. 
 If
 $A$ belongs to $ [\Naturals]^{\aleph_0}$, but not necessarily to ${\mathfrak
   B}_\theta$, 
and $G$ is $\Qposet({\mathfrak S})$ generic over $V$
then there are infinitely many integers $j\in A$ such that
$F_{\mathfrak S}[G](j)  = F_{\theta}[H](j)$.
\end{lemma}
\begin{proof}
  This is a standard use of genericity.
\end{proof}
\begin{theor}\label{t:mmm}
Given any regular, uncountable cardinal $\kappa $, it is consistent relative to the consistency of set theory that
Martin's Axiom holds, $2^{\aleph_0} = \kappa$ and there is a nowhere
trivial automorphism of ${\mathcal P}(\Naturals)$ modulo the finite sets.
\end{theor}
\begin{proof}
  Let $C\subseteq \kappa$ be a closed unbounded set such that $\kappa\setminus
  C$ is unbounded. Construct a finite support iteration
  $\{\Poset_\xi\}_{\xi\in C}$ so that Martin's Axiom is forced to hold by
  the iteration on $\kappa\setminus C$ and such that a tower of permutations $
  \{(A_\xi,F_\xi, {\mathfrak B}_\xi)\}_{\xi\in \kappa}$ is generically
  constructed along $C$ such that ${\mathfrak B}_{\xi} = \pomega\cap
  V[\Poset_\xi\cap G]$ where $G$ is $\Poset_\kappa$ generic over $V$.
Then let $\Psi$ be the automorphism of $\pomegaf$ defined by $\Psi([A]) = [B]$
  if and only if there is some $\xi \in \kappa$ such that $F_\xi(A) \equiv^*
  B$ and $A \in {\mathfrak B}_\xi$. To see that $\Psi$ is nowhere trivial
  suppose that $\Psi$ is induced by $\psi$ on $A$. Let $\xi\in C$ be an ordinal
  large enough that $A$ and $\psi$ both belong to $V[G\cap \Poset_\xi]$. 
Using Lemma~\ref{l:kl} it is possible to find in $V[G\cap \Poset_{\xi + 1}]$
  an infinite set $Z$ such that $\psi(Z) \cap F_\xi(Z) = \emptyset$. However,
  there is no guarantee that $Z$ belongs to ${\mathfrak B}_\xi$. Let $\rho$ be
  the first member of $C$ greater than $\xi$. Then 
 $Z$ does belong to ${\mathfrak B}_\rho$. From Lemma~\ref{l:kl+} it follows
  that $F_\rho(j) = F_\xi(j)$ for infinitely many $j$ belonging to $Z$.
Therefore $F_\rho(Z) \setminus \psi(Z)$ is infinite contradicting that $\psi$
  induces $\Psi$ on $A$ and $[F_\rho(Z)] = \Psi([Z])$.
\end{proof}

It should be noted that Theorem~\ref{t:mmm} would be of interest even if
Martin's  Axiom did not hold in the model constructed since it  would still
 provide
a method for constructing models of set theory with nowhere trivial
automorphisms of 
$\pomegaf$ and the continuum arbitrarily large. 

\section{Ruining automorphisms with Silver reals}
  \begin{defin}\label{d:ss}
 Suppose that $\Rationals$ and $\Poset$ are partial orders such that
 $\Rationals$ is completely embedded in $\Poset$ and that $\Rationals$ is
 Suslin (see \cite{ba.ju.book}).  Then $\Poset$ 
will be said to be {\em sufficiently
 Suslin} over $\Rationals$ if for every  $\Poset$-name $g$ for a
 function from $\omega$ to $2$
  there is a dense set of 
$ p \in \Poset$ such that,
$$\left\{(q,f) \in \Rationals \times 2^\omega | (\exists p'\in \Poset)
p' \geq q \AND
p' \geq p \AND p'\forces{\Poset}{g = \check{f}}\right\}$$ 
is analytic. 
\end{defin}
  \begin{lemma}\label{l:sp}
Any countable support iteration of a combination of Silver reals and Sacks
reals is sufficiently Suslin over its first coordinate.    
  \end{lemma}
\begin{proof}
Let the iteration be obtained from the sequence $\{\Poset_\xi\}_{\xi\in \eta}$
where each successor stage is constructed by using one of the mentioned
partial orders and let the  $\Poset_\eta$-name $g$  be given.
Let ${\mathfrak M} \prec (H(\card{\Poset_\eta}^+), \in)$ be a countable
elementary submodel containing $\Poset_\eta$ and  $g$. 
Let ${\mathfrak M}\cap \eta = \Gamma$, let $\{\gamma_i\}_{i\in\omega}$
enumerate $\Gamma$ and let $\Gamma_j =\{\gamma_i\}_{i\in j} $.
Standard fusion arguments will allow the construction of a  family
$\{p_\tau\}_{\tau : \Gamma_n \to 2^{<n}} $ 
 such that:
\begin{enumerate}
\item If $\tau(\gamma) \subseteq \tau'(\gamma)$ for each $\gamma$  then $p_\tau
  \leq p_{\tau'}$. 
\item If $\tau : \Gamma_n \to 2^{<n}$ then $p_\tau$ decides
the value of $g(n)$.
\item If $p(\gamma)$ is defined to be $$\bigwedge_{n\in\omega}\left(
  \bigvee_{\tau:\Gamma_n \to 2^{<n}\AND p_\tau\restriction \gamma \in
    G}p_\tau\right)$$ 
  then $p \in \Poset$.
\end{enumerate}
It follows that, for each $q \in \Poset_0$ there is some 
$p' \geq p $ such that $p'(0) \geq q$ and
$p'\forces{\Poset}{g = \check{f}}$ if and only if, letting 
$S(f,n)  = \{\tau | \tau:\Gamma_n \to 2^{<n}\AND p_\tau\forces{\Poset}{g(n) =
        \check{f}(\check{n})}\AND p_\tau(0) \geq q\}$, the set of conditions
$\{
  \bigvee_{\tau \in S(f,n)}
p_\tau\}_{n\in\omega}$ 
has a proper lower bound. For Silver and Sacks forcings  checking for
a lower bound is easily seen to be $\Sigma^1_1$ in the parameter defining the
fusion sequence; in fact, it is Borel.
\end{proof}
\begin{lemma}\label{l:case1}
If $\Poset$ is as in lemma~\ref{l:sp} and
 $\Psi$ is an automorphism of $\pomegaf$ which is not
trivial on any member of
$ [\omega]^{\aleph_0}$ 
then it is not possible to extend $\Psi$ to an automorphism of $\pomegaf$ in
 any generic extension by  $\Sacks*\Poset$.
\end{lemma}  
\begin{proof}
Assuming the lemma is false, it
is possible to find a condition $(s,p) \in \Sacks * \Poset$ and an
$\Sacks * \Poset$-name for a set of integers $Z$ such that
$$(s,p)\forces{\Sacks * \Poset}{\Psi^*\text{ is a lifting of an extension of
    }\Psi\text{ and }\Psi^*(X_G) = Z}$$
where $X_G$ is a name for the generic subset of $\omega$ added by $\Sacks$.
Let $g= \chi_Z$ be the characteristic function of $Z$.  Using the fact that
$\Sacks*\Poset$ is sufficiently Suslin, find $(s',p')\in \Sacks * \Poset$ such
that $(s',p') \geq (s,p)$ and 
$$\left\{(\bar{s},f) \in \Sacks \times 2^\omega | (\exists q \in \Sacks *
  \Poset)q(0) \supseteq \bar{s}  \AND
q \geq (s',p') \AND q\forces{\Sacks *\Poset}{g = \check{f}}\right\}$$ 
is analytic.
Let $A$ be an infinite set of integers disjoint from the domain of $s'$ such
that $\dom(s') \cup A$ has infinite complement in the integers.
It follows that if $R$ is defined to be
$\left\{(t,f) \in 2^A\times 2^\omega : (\exists q \geq (t\cup s',p'))
q\forces{\Sacks * \Poset}{g=\check{f}}\right\}$ then $R$ is
analytic. 

There are two possibilities. First, suppose that
the domain of $R$ is all of $2^A$. In this case it is
possible to find a continuous function $S$ defined on a comeagre subset of 
${\mathcal P}(A)$ such that $R(t,S(t^{-1}\{1\}))$ holds for all $t$ in the
domain of $R$. 
Now let $\psi: \dom(S) \to {\mathcal P}(\Psi^*(A))$ be defined by
$$\psi(W) = \{n \in \Psi^*(A) : S(W)(n)= 1\}$$
and observe that $\psi$ is continuous. Furthermore, 
 $\psi(W)  \equiv^* 
\Psi^*(W)$. To see this,
let $q \in \Sacks *\Poset$ be any condition witnessing  
that $R(\chi_W\restriction A, S(W))$ holds. 
Observe that not only does
$q$ force $X_G\cap A = W$ but also
$q\forces{\Sacks * \Poset}{Z \cap \Psi^*(A) = \psi(W)}$.
Hence,
$q\forces{\Sacks * \Poset}{\psi(W)
  =\Psi^*(X_G)\cap \Psi^*(A) \equiv^* \Psi^*(X_G\cap A) \equiv^*  
\Psi^*(W)}$. 
Therefore, since $\psi(W)$ and $\Psi^*(W)$ are sets in the ground
model, it follows from the absoluteness of $\equiv^*$
 that    $\psi(W)\equiv^* 
\Psi^*(W)$.
 But now it follows that $\Psi$ is trivial on $A$ by
Lemma~1 of \cite{veli.def}.  

In the other case there is some $t: A \to 2$ such that there is no $f$ such
that $R(t,f)$ holds. 
This implies that for every $q \geq (s'\cup t,p')$ there are
infinitely many integers in $\Psi^*(A)$ whose membership in $Z$ is not decided
by $q$. Genericity implies that
$(s'\cup t,p')\forces{\Sacks * \Poset}{Z\cap \Psi^*(A) \not\equiv^*
\Phi^*(t^{-1}\{1\})}$ which 
is a contradiction to the assumption that
$(s,p)\forces{\Sacks * \Poset}{\Psi^*(X_G) = Z}$.
\end{proof}
\begin{theor}
Let $V$ be a model of \CH \. 
If $\Poset_{\omega_2}$ is the countable support iteration of 
partial orders as in Lemma~\ref{l:sp} such that $\Poset_{\alpha+1} =
\Poset_{\alpha}*\Sacks$ for a stationary
set of $\alpha \in \omega_2$ then $V[G]$ is a model where every
automorphism of $\pomegaf$ is somewhere trivial for every $G\subseteq
\Poset_{\omega_2}$ which is generic over $V$. (Hence, it is consistent that
 every
automorphism of $\pomegaf$ is somewhere trivial and ${\mathfrak d} =\aleph_1$.) 
\end{theor}
\begin{proof}
Any automorphism can be reflected 
on a closed unbounded subset of $\omega_2$ consisting of ordinal of
uncountable cofinality. If $\alpha$ is any such ordinal
such that $\Poset_{\alpha+1} =
\Poset_{\alpha}*\Sacks$ then  Lemma~\ref{l:case1} can be
applied to show that the automorphism can not be extended any further.
\end{proof}

Finally, it should be remarked that the hypothesis of Lemma~\ref{l:sp} can be
extended to include various other partial orders. However, in light of the
lack of immediate applications and the technical
difficulties required to establish this, 
the proof will be provided elsewhere.


\end{document}